\documentclass[12pt]{article}

      \usepackage{latexsym}
         \usepackage[reqno, namelimits, sumlimits]{amsmath}
         \usepackage{amssymb, amsfonts}
         \usepackage{amsthm}
 \newtheorem{theorem}{Theorem}[section]

\begin{document}

\begin{center} {\Large\bf  A certain necessary condition of potential  blow up for Navier-Stokes equations}\\
  \vspace{1cm}
 {\large
  G. Seregin}\footnote{OxPDE, Mathematical Institute, University of Oxford, UK}

 \end{center}

 \vspace{1cm}

 \vspace{1cm}
 \noindent
 {\bf Abstract }  %We show that the spatial $L_3$-norm of the velocity field tends to $\infty$ as time $t$ approaches $T$ provided $T$ is a potential blow up time.
We show that a necessary condition for $T$ to be a potential blow up time is
$\lim\limits_{t\uparrow T}\|v(\cdot,t)\|_{L_3}=\infty$.

 \vspace {1cm}

\noindent {\bf 1991 Mathematical subject classification (Amer.
Math. Soc.)}: 35K, 76D.


\setcounter{equation}{0}
\section{Main Result}

%In the present paper, we address the following question.
Consider the Cauchy problem
for the classical Navier-Stokes system
\begin{equation}\label{m1}
    \partial_t v +v \cdot \nabla v -\Delta v=-\nabla q,\qquad {\rm div} v=0,
\end{equation}
describing the flow of a viscose incompressible fluid in $Q_+=\mathbb R^3\times ]0,\infty[$ with the initial condition
\begin{equation}\label{m2}
    v|_{t=0}=a
\end{equation}
in $\mathbb R^3$. Here, as usual, $v$ and $q$ stand for the velocity field and for the pressure field, respectively. For simplicity, let us assume
\begin{equation}\label{m3}
    a\in C^\infty_{0,0}(\mathbb R^3)\equiv \{v\in C^\infty_0(\mathbb R^3):\,\,{\rm div} v=0\}.
\end{equation}

Many years ago, in 1934, J. Leray showed in his celebrated paper \cite{Le} that  problem (\ref{m1})--(\ref{m3}) has at least one weak solution obeying the  global energy inequality
\begin{equation}\label{m4}
    \frac 12\int\limits_{\mathbb R^3}|v(x,t)|dx+\int\limits^t_0\int\limits_{\mathbb R^3}|\nabla v|^2dxdt'\leq \frac 12 \int\limits_{\mathbb R^3}|a|^2dx
\end{equation}
for all positive $t$. This solution is smooth and unique for sufficiently small values of $t$. The first instant of time $T$ when singularities occur is called a {\it blow up} time. By definition, $z_0=(x_0,t_0)$ is a \textit{regular} point of $v$ if it is essentially bounded in a nonempty parabolic ball with the center at the point $z_0$ \footnote{$Q(z_0,r)=B(x_0,r)\times ]t_0-r^2,t_0[$ is a parabolic ball of radius $r$ centered at the point $z_0$.} The point $z_0$ is  \textit{singular} if it is not  regular.

To the best of our knowledge, it is unknown whether there exists
an energy solution to the Cauchy problem (\ref{m1})--(\ref{m3}) with a finite time blow up.

However, J. Leray proved certain necessary conditions for $T$ to be a blow up time. They are as follows. Assume that $T$ is a  blow up time, then, as it has been shown in the above mentioned  paper,
for any $3<m\leq \infty$, there exists a constant $c_m$, depending on $m$ only, such that
\begin{equation}\label{m5}
    \|v(\cdot,t)\|_m\equiv \|v(\cdot,t)\|_{m,\mathbb R^3}\equiv\Big(\int\limits_{\mathbb R^3}|v(x,t)|^mdx\Big)^\frac 1m\geq \frac {c_m}{(T-t)^\frac {m-3}{2m}}
\end{equation}
for all $0<t<T$.

Here, we address the critical case $m=3$, for which
%is of great interest since the norm $\|\cdot\|_3$ is critical in the sense that it is invariant with respect to the Navier-Stokes scaling:
%\begin{equation}\label{m11}
%    v^\lambda(x,t)=\lambda v(\lambda x,\lambda^2t),\qquad q^\lambda(x,t)=\lambda^2 q(\lambda %x,\lambda^2t)
%\end{equation}
%with a positive parameter $\lambda$, there has been proved
a weaker statement
\begin{equation}\label{m6}
    \limsup\limits_{t\to T-0}\|v(\cdot,t)\|_3=\infty
\end{equation}
has been proven in \cite{ESS4}. The aim of this paper and several previous papers  of the author is to improve (\ref{m6}). At the moment, the best improvement of (\ref{m6})  is given by the following theorem.

\begin{theorem}\label{mainresult} Let $v$ be an energy solution to the Cauchy problem
(\ref{m1}) and (\ref{m2}) with the  initial data satisfying (\ref{m3}). Let $T>0$ be a finite blow up time. Then \begin{equation}\label{m7}
  \lim\limits_{t\to T-0}\|v(\cdot,t)\|_3=\infty
\end{equation} holds true.
\end{theorem}

Now, let us briefly outline a proof that relays  upon ideas developed in  %previous papers of the author on this matter, see
\cite{S5}-\cite{S7}. In particular, in \cite{S5}, a certain type of scaling has been invented, which, after passing to the limit, gives a special non-trivial solution to the Navier-Stokes equations provided %under the assumption
 there is a finite time blow up. In \cite{S6} and \cite{S7}, it has been shown  that the same type of scaling and blowing-up can produce  the so-called Lemarie-Rieusset local energy solutions. It turns out to be that %The great advantage of such kind of solutions is that %they possess a spatial decay, which is one of two %main conditions for %applications of
the backward uniqueness technique is still   applicable to them.
%The known way is to reduce the problem either to Liouville type theorems for bounded ancient solutions,  see \cite{KNSS} and \cite{SS3}, or to backward uniqueness for the heat operator with lower order term, see \cite{ESS4}. So far, the experience  shows  us that
%working with scale-invariant norms it is preferable to utilize the second approach.
Although the theory of backward uniqueness itself is relatively well understood, its realization is not an easy task and based on delicate  regularity results for the Navier-Stokes equations. Actually, there are two main points to verify: solutions, produced by scaling and blowing-up, vanish at the last moment of time and have a spatial decay. The first property is easy when working with $L_3$-norm while the second one is harder. However, under certain restrictions,
the required decay is a consequence of the Lemarie-Rieusset theory. So, the main technical part of the whole procedure is to show that scaling and blowing-up lead to local energy solutions.
On that way, a lack of compactness of initial data of scaled solutions %are compact
in $L_{2,\rm{loc}}$ is the main obstruction.  %in showing %It is not so difficult to show
%that the limit solution is a local energy solution. % is a of %if one can show that
%initial data of scaled solutions %are compact
%in $L_{2,\rm{loc}}$.
This is why the same theorem for a stronger scale-invariant norm of the space $H^\frac 12$ is easier. The reason for that is a compactness of  the corresponding embedding, see  \cite{RS} and \cite{S6}.

In  this paper, we are going  to show that, despite of a lack of compactness in $L_3$-case, the limit of the sequence  of scaled solutions is still a local energy solution, for which a spatial decay takes place.
%how
%one  can  produce a local energy solution as the limit of the sequence  of scaled solutions the %limit of the sequence  of scaled solutionseven without above mentioned compactness.
Technically, this can be done by splitting each scaled solution into two parts. The first one is a solution to a non-linear problem but with zero initial data while the second one is a solution of a linear problem with weakly converging nonhomogeneous initial data.

\setcounter{equation}{0}
\section{Estimates of Scaled Solutions}

%There is a common part when proving  Theorem \ref{mainresult} and Proposition \ref{l3} and it is as follows.
Assume that  our statement is false and there exists an increasing sequence $t_k$ converging to $T$ as $k\to \infty$ such that
\begin{equation}\label{sb1}
    \sup\limits_{k\in \mathbb N}\|v(\cdot,t_k)\|_{3}=M<\infty.
\end{equation}
%where $\mathcal B$ is either $H^\frac 12$ or $L_3$.

 By the definition of a blow up time for energy solutions, %For the case of energy solutions, it is known that   if  $T$ is a blow up time, then
 there exists at least one singular point at time $T$. Without loss of generality, we may assume that it is $(0,T)$. Moreover, the blow-up profile has the finite $L_3$-norm, i.e.,

\begin{equation}\label{sb2}
    \|v(\cdot,T)\|_{3}<\infty.
\end{equation}

Let us scale $v$ and $q$ %in the following way
so that
\begin{equation}\label{sb3}
    u^{(k)}(y,s)=\lambda_kv(x,t),\quad p^{(k)}(y,s)
    =\lambda_k^2q(x,t),\quad (y,s)\in \mathbb R^3\times ]-\lambda_k^{-2}T,0[,
\end{equation}
where
$$x=\lambda_ky,\qquad t=T+\lambda_k^2s,$$
$$\lambda_k=\sqrt{\frac {T-t_k}{S}}$$
and a positive parameter $S<10$ will be defined later.

By the scale invariance of  $L_3$-norm, $u^{(k)}(\cdot,-S)$ is uniformly bounded in $L_3(\mathbb R^3)$, i.e.,
% $scale-invariant, scaled functions have the following property
\begin{equation}\label{sb4}
  \sup\limits_{k\in \mathbb N}  \|u^{(k)}(\cdot,-S)\|_{3}=M<\infty.
\end{equation}

Let us decompose our scaled solution $u^{(k)}$ %and $p^{(k)}$
into two parts: $u^{(k)}=v^{(k)}+ w^{(k)}$.
Here, $w^{(k)}$ is a solution to the Cauchy  problem for the Stokes system:
\begin{eqnarray}\label{f1}
% \nonumber to remove numbering (before each equation)
   \partial_t w^{(k)}-\Delta w^{(k)}=-\nabla r^{(k)},\quad {\rm div}\, w^{(k)}=0 \quad{\rm in} \quad \mathbb R^3\times ]-S,0[,
     \nonumber\\
  w^{(k)}(\cdot,-S)=u^{(k)}(\cdot,-S).
\end{eqnarray}
Apparently, (\ref{f1})
%It well-known that this problem problem
can be reduced to the Cauchy problem for the heat equation
so that the pressure $ r^{(k)}=0$ and  $w^{(k)}$ can be worked out with the help of %expressed via
the heat potential. The %following
 estimate below is well-known, see, for example \cite{Kato},
\begin{equation}\label{f2}
  \sup\limits_k  \{\|w^{(k)}\|_{L_5(\mathbb R^3\times ]-S,0[}+\|w^{(k)}\|_{L_{3,\infty}(\mathbb R^3\times ]-S,0[}\}\leq c(M)<\infty.
\end{equation}
It is worthy to note that, by the scale invariance, $c(M)$ in  (\ref{f2}) is independent of $S$.

 As to $v^{(k)}$, it is a solution to the Cauchy problem for the following perturbed Navier-Stokes system
\begin{eqnarray}\label{f3}
% \nonumber to remove numbering (before each equation)
   \partial_t v^{(k)}+{\rm div}(v^{(k)}+w^{(k)})\otimes(v^{(k)}+w^{(k)}) -\Delta v^{(k)}=-\nabla p^{(k)},\nonumber\\ {\rm div}\, v^{(k)}=0 \quad{\rm in} \quad \mathbb R^3\times ]-S,0[,
     \\
  v^{(k)}(\cdot,-S)=0\nonumber .
\end{eqnarray}

Now, our aim is to show that, for a suitable choice of $-S$, we can prove unform estimates of $v^{(k)}$ and $p^{(k)}$ in certain spaces, pass to the limit as $k\to \infty$, and conclude that the limit functions $u$ and $p$ are a local energy solution to the Cauchy problem for the Navier-Stokes system in $\mathbb R^3\times ]-S,0[$ associated the initial data, generated by  the weak $L_3$-limit of the sequence $u^{(k)}(\cdot,-S)$.
%  in the sense of the definition accepted in \cite{KS}, see \cite{LR1} for an earlier version.

 Let us start with estimates of solution to (\ref{f3}). First of all, we know the formula for the  pressure:
\begin{equation}\label{f4}
  p^{(k)}(x,t)=-\frac 13  |u^{(k)}(x,t)|^2 +\frac 1{4\pi}\int\limits_{\mathbb R^3}K(x-y):u^{(k)}(y,t)\otimes u^{(k)}(y,t)dy,
\end{equation}
where
%$u^{(k)}=v^{(k)}+ w^{(k)}$ and and
$K(x)=\nabla^2(1/|x|)$.

Next, we may decompose the pressure in the same way as it has been done in \cite{KS}. For $x_0\in \mathbb R^3$ and for $x\in B(x_0,3/2)$, we let
\begin{equation}\label{f5}
   {p} ^{(k)}_{x_0}(x,t)\equiv p^{(k)}(x,t)-c^{(k)}_{x_0}(t)=p^{1(k)}_{x_0}(x,t)+p^{2(k)}_{x_0}(x,t),
\end{equation}
where
$$p^{1(k)}_{x_0}(x,t)=-\frac 13|u^{(k)}(x,t)|^2+\frac 1{4\pi}\int\limits_{B(x_0,2)}K(x-y): u^{(k)}(y,t)\otimes u^{(k)}(y,t)dy,$$
%$$q^{2(k)}_{x_0,r,R}(x,t)=\frac 1{4\pi}\int\limits_{B(x_0,2R)\setminus B(x_0,r)}(K(x-y)-K(x_0-y)):u^{(k)}(y,t)\otimes u^{(k)}(y,t)dy,$$
$$p^{2(k)}_{x_0}(x,t)=\frac 1{4\pi}\int\limits_{\mathbb R^3\setminus B(x_0,2)}(K(x-y)-K(x_0-y)):u^{(k)}(y,t)\otimes u^{(k)}(y,t)dy,$$
$$c^{(k)}_{x_0}(t)=\frac 1{4\pi}\int\limits_{\mathbb R^3\setminus B(x_0,2)}K(x_0-y):u^{(k)}(y,t)\otimes u^{(k)}(y,t)dy.$$

Using the similar arguments as in \cite{LR1}, one can derive estimates for the above counterparts of the pressure. Here, they are:
\begin{equation}\label{f6}
\|p^{1(k)}_{x_0}(\cdot,t)\|_{L_\frac 32(B(x_0,3/2))}\leq c(M)(   \|v^{(k)}(\cdot,t)\|^2_{L_ 3(B(x_0,2))}+1),
\end{equation}

\begin{equation}\label{f8}
 \sup\limits_{B(x_0,3/2)}| p^{2(k)}_{x_0}(x,t) |\leq c(M)(   \|v^{(k)}(\cdot,t)\|^2_{L_{ 2,{\rm unif}}}+1),
\end{equation}
where
$$\|g\|_{L_{ 2,{\rm unif}}}=\sup\limits_{x_0\in \mathbb R^3}\|g\|_{L_{ 2}(B(x_0,1))}.$$

We further let
$$ \alpha(s)=\alpha(s;k,S)=\|v^{(k)}(\cdot,s)\|^2_{2,{\rm unif}},$$$$ \beta(s)=\beta(s;k,S)=\sup\limits_{x\in R^3}\int\limits_{-S}^s\int\limits_{B(x,1)}|\nabla v^{(k)}|^2dyd\tau.$$

From (\ref{f6}), (\ref{f8}), we find the estimate of the scaled pressure
\begin{equation}\label{f9}
\delta(0) \leq c(M)\Big[\gamma(0)+\int\limits_{-S}^0(1+\alpha^\frac 32(s))ds\Big],
\end{equation}
with some positive constant $c(M)$ independent of $k$ and $S$. Here, $\gamma$ and $\delta$ are defined as
$$\gamma(s)=\gamma(s;k,S)=\sup\limits_{x\in R^3}\int\limits_{-S}^s\int\limits_{B(x,1)}| v^{(k)}(y,\tau)|^3dyd\tau$$ and $$\delta(s)=\delta(s;k,S)=\sup\limits_{x\in R^3}   \int\limits_{-S}^s\int\limits_{B(x,3/2)}|p^{(k)}(y,\tau)-c^{(k)}_x(\tau)|^\frac 32dy\,d\tau,
$$
%with some function $c^{(k)}_x\in L_\frac 32(-S,0)$ and
 respectively. It is known that an upper bound for $\gamma$ can be given by the known multiplicative inequality
\begin{equation}\label{f10}
  \gamma(s)\leq c\Big(\int\limits_{-S}^s\alpha^3(\tau)d\tau\Big)^\frac 14\Big(\beta(s)+  \int\limits_{-S}^s\alpha(\tau)d\tau\Big)^\frac 34.
\end{equation}

Fix $x_0\in \mathbb R^3$ and a smooth non-negative function $\varphi$ such that
$$\varphi=1\quad {\rm in }\quad B(1),\qquad {\rm spt}\,\subset B(3/2)$$
and let $\varphi_{x_0}(x)=\varphi(x-x_0)$.

Since the function $v^{(k)}$ is smooth on $[-S,0[$, all our further actions are going to be legal. In particular, we may write down the following energy identity
$$\int\limits_{\mathbb R^3 }\varphi^2_{x_0}(x)|v^{(k)}(x,s)|^2dx+
2\int\limits^s_{-S}\int\limits_{\mathbb R^3 }\varphi^2_{x_0}|\nabla v^{(k)}|^2dxd\tau=$$
 $$=\int\limits^s_{-S}\int\limits_{\mathbb R^3 }\Big[| v^{(k)}|^2\Delta\varphi^2_{x_0}+
  v^{(k)}\cdot \nabla \varphi^2_{x_0}(| v^{(k)}|^2+2p^{(k)}_{x_0})\Big]
 dxd\tau+$$
 $$+\int\limits^s_{-S}\int\limits_{\mathbb R^3 }\Big[ w^{(k)}\cdot \nabla \varphi^2_{x_0}| v^{(k)}|^2
 +2\varphi^2_{x_0}w^{(k)}\otimes(w^{(k)}+v^{(k)}):\nabla v^{(k)}+$$
 $$+2w^{(k)}\cdot v^{(k)}(w^{(k)}+v^{(k)})\cdot \nabla \varphi^2_{x_0} \Big]dx d\tau =I_1+I_2.$$

 The first term $I_1$ is estimated  with the help of the H\"older inequality, multiplicative inequality (\ref{f10}), and bounds (\ref{f6}), (\ref{f8}). So, we find
 $$I_1\leq c(M)\Big[\int\limits^s_{-S}(1+\alpha(\tau)+\alpha^\frac 32(\tau))d\tau+$$
 $$+\Big(\int\limits_{-S}^s\alpha^3(\tau)d\tau\Big)^\frac 14\Big(\beta(s)+  \int\limits_{-S}^s\alpha(\tau)d\tau\Big)^\frac 34\Big].$$
 Now, let us evaluate the second term
 $$I_2\leq c\int\limits_{-S}^s\|v^{(k)}(\cdot,\tau)\|^2_{L_3(B(x_0,3/2))}
 \|w^{(k)}(\cdot,\tau)\|_{L_3(B(x_0,3/2))}d\tau+$$
 $$+c\int\limits_{-S}^s\Big(\int\limits_{B(x_0,3/2)}|w^{(k)}|^5dx\Big)^\frac 15\Big(\int\limits_{B(x_0,3/2)}|v^{(k)}|^\frac 54|\nabla v^{(k)}|^\frac 54dx\Big)^\frac 45d\tau+$$
 $$+c\beta^\frac 12(s)\Big(\int\limits_{-S}^s\int\limits_{B(x_0,3/2)}|w^{(k)}|^4dxd\tau\Big)^\frac 12d\tau+$$
 $$+c\int\limits_{-S}^s\|v^{(k)}(\cdot,\tau)\|_{L_3(B(x_0,3/2))}
 \|w^{(k)}(\cdot,\tau)\|^2_{L_3(B(x_0,3/2))}d\tau.$$
Taking into account  and applying H\"older inequality several times(\ref{f2}), we find
$$I_2\leq c(M)\gamma^\frac 23(s)(s+S)^\frac 13+$$$$+c \int\limits_{-S}^s\Big(\int\limits_{B(x_0,3/2)}|w^{(k)}|^5dx\Big)^\frac 15\Big(\int\limits_{B(x_0,3/2)}|\nabla v^{(k)}|^2dx\Big)^\frac 12\times$$$$\times \Big(\int\limits_{B(x_0,3/2)}| v^{(k)}|^\frac {10}3dx\Big)^\frac 3{10}d\tau
%\Big(\int\limits_{B(x_0,3/2)}|w^{(k)}|^5dxd\tau\Big)^\frac 15$$
+c(M)\beta^\frac 12(s)(s+S)^\frac 1{10}+$$$$+c(M)\gamma^\frac 13(s)(s+S)^\frac 23.$$
It remains to use another known multiplicative inequality
$$\Big(\int\limits_{B(x_0,3/2)}|v^{(k)}(x,s)|^\frac{10}3dx\Big)^\frac 3{10}\leq c\Big(\int\limits_{B(x_0,3/2)}|v^{(k)}(x,s)|^2dx\Big)^\frac 15\times$$
$$\times \Big(\int\limits_{B(x_0,3/2)}(|\nabla v^{(k)}(x,s)|^2+|v^{(k)}(x,s)|^2dx\Big)^\frac 3{10}$$
and to conclude that %Then we have
$$I_2\leq c(M)\gamma^\frac 23(s)(s+S)^\frac 13+c(M)\beta^\frac 12(s)(s+S)^\frac 1{10}+c(M)\gamma^\frac 13(s)(s+S)^\frac 23+$$
$$+c\Big(\beta(s)+\int\limits_{-S}^s\alpha(\tau)d\tau)\Big)^\frac 45
\times\Big(\int\limits_{-S}^s\alpha(\tau)\|w^{(k)}(\cdot,\tau)\|^5_{L_{5,{\rm unif}}}d\tau\Big)^\frac 15.$$
Finally, we find
$$ \alpha(s)+\beta(s)\leq c(M)\Big[(s+S)^\frac 15+$$
\begin{equation}\label{f11}
  + \int\limits^s_{-S}\Big(\alpha(\tau)(1+\|w^{(k)}(\cdot,\tau)\|^5_{L_{5,{\rm unif}}})+\alpha^3(\tau)\Big)\,d\tau\Big],
 \end{equation}
which is valid for any $ s\in [-S,0[$ and for some positive constant $c(M)$ independent of $k$, $s$, and $S$.

It is not so difficult to show that there is a  positive constant  $S(M)$ such that
\begin{equation}\label{f12}
    \alpha(s)\leq \frac 1{10}
\end{equation}
for any $s\in ]-S(M),0[$. In turn, the latter will also imply that
\begin{equation}\label{f13}
   \alpha(s)\leq c(M)(s+S)^\frac 15
\end{equation}
for any $s\in ]-S(M),0[$.

To see how this can be worked out, let us assume
\begin{equation}\label{f14}
    \alpha(s)\leq 1
\end{equation}
for $-S\leq s< s_0\leq 0$. Then  (\ref{f11}) yields
\begin{equation}\label{f15}
    \alpha(s)\leq c(M)((s+S)^\frac 15+y(s))
\end{equation}
for the same $s$. Here,
$$y(s)=\int\limits^s_{-S}\alpha(\tau)(2+g(\tau))d\tau, \qquad g(s)=\|w^{(k)}(\cdot,s)\|^2_{L_5(\mathbb R^3)}.$$
The function $y(s)$ obeys the differential inequality
\begin{equation}\label{f16}
    y'(s)\leq c(M)(2+g(s))((s+S)^\frac 15+y(s))
\end{equation}
for $-S\leq s< s_0\leq 0$. After integrating (\ref{f16}), we find
\begin{equation}\label{f17}
    y(s)\leq c(M)\int\limits^s_{-S}\Big((\tau+S)^\frac 15(2+g(\tau))\exp{\Big\{c(M)\int\limits^s_\tau(2+g(\vartheta))\Big\}}d\vartheta\Big)d\tau
\end{equation}
for $-S\leq s< s_0\leq 0$. Taking into account estimate (\ref{f2}), we derive from (\ref{f17}) the following bound
\begin{equation}\label{f18}
    y(s)\leq c_1(M)(s+S)^\frac 15
\end{equation}
for $-S\leq s< s_0\leq 0$ and thus
\begin{equation}\label{f19}
    y(s)\leq c_1(M)S^\frac 15
\end{equation}
for the same $s$.

Now, let us pick up $S(M)>0$ so small that
\begin{equation}\label{f20}
    c(M)(1+c_1(M))S^\frac 15(M)=\frac 1{20}.
\end{equation}
We claim that, for such a choice of S(M), statement (\ref{f12}) holds true.
Indeed, assume that it is false. Then since $\alpha(s)$ is a continuous function
on $[-S,0[$ and $\alpha(0)=0$, there exists $s_0\in ]-S,0[$ such that
$0\leq\alpha(s)<\frac 1{10}$ for all $s\in ]-S,s_0[$ and $\alpha(s_0)=\frac 1{10}$. In this case, we may use first (\ref{f19}) and then (\ref{f15}), (\ref{f20}) to get
$$\alpha(s)\leq c(M)(1+c_1(M))S^\frac 15(M)=\frac 1{20}$$
for $s\in ]-S,s_0[$. This leads to a contradiction and, hence, (\ref{f12}) has been proven.
It remains to use (\ref{f15}) and (\ref{f18}) with $s_0=0$ in order to establish (\ref{f13}).

\setcounter{equation}{0}
\section{Limiting Procedure}

As to $w^{(k)}$, it is defined by the solution formula
$$w^{(k)}(x,t)=\frac 1{(4\pi(s+S))^\frac 32}\int\limits_{\mathbb R^3}\exp{\Big(-\frac {|x-y|^2}{4(s+S)}\Big)}u^{(k)}(y,-S)dy.$$
Moreover, by standard localization arguments, the following estimate can be derived:
$$\sup\limits_{-S<s<0}\sup\limits_{x_0\in\mathbb R^3}\|w^{(k)}(\cdot,s)\|^2_{L_2(B(x_0,1))}+
$$$$+\sup\limits_{x_0\in\mathbb R^3}\int\limits_{-S}^0\int\limits_{B(x_0,1)}|\nabla w^{(k)}(y,s)|^2dyds\leq c(M)<\infty.
$$

Obviously, $w^{(k)}$ and all its derivatives converge to $w$ and to its corresponding derivatives uniformly in  sets of the form
$\overline{B}(R)\times [\delta,0]$ for any $R>0$ and for any $\delta\in ]-S,0[$. The limit function satisfies the same representation formula
$$w(x,t)=\frac 1{(4\pi(s+S))^\frac 32}\int\limits_{\mathbb R^3}\exp{\Big(-\frac {|x-y|^2}{4(s+S)}\Big)}a_0(y)dy,$$
in which $a_0$ is the weak $L_3(\mathbb R^3)$-limit of the sequence $u^{(k)}(\cdot,-S)$. The function $w$ satisfies the uniform local energy estimate
$$\sup\limits_{-S<s<0}\sup\limits_{x_0\in\mathbb R^3}\|w(\cdot,s)\|^2_{L_2(B(x_0,1))}+
$$$$+\sup\limits_{x_0\in\mathbb R^3}\int\limits_{-S}^0\int\limits_{B(x_0,1)}|\nabla w(y,s)|^2dyds\leq c(M)<\infty.
$$
The important fact, coming from the solution formula, is as follows:
\begin{equation}\label{l1}
    w\in C([-S,0];L_3(\mathbb R^3))\cap L_5(\mathbb R^3\times ]-S,0[).
\end{equation}

Next,  the uniform  local energy estimate for the sequence $u^{(k)}$ (with respect to $k$) can be deduced from the estimates above.  This allows us to exploit %is uniformly bounded  and we can use
the limiting procedure explained in \cite{KS} in details. As a result, one can selected a subsequence, still denoted by $u^{(k)}$,
with the following properties:

 for any $a>0$,
\begin{equation}\label{lp1}
   u^{(k)}\to u
\end{equation}
weakly-star in $L_\infty(-S,0;L_2(B(a)))$ and strongly in $L_3(B(a)\times ]-S,0[)$ and in
$C([\tau,0];L_\frac 98(B(a)))$ for any $-S<\tau<0$;
\begin{equation}\label{lp2}
    \nabla u^{(k)}\to \nabla u
\end{equation}
weakly in $L_2(B(a)\times ]-S,0[)$;
\begin{equation}\label{lp3}
    t\mapsto\int\limits_{B(a)}u^{(k)}(x,t)\cdot w(x) dx\to t\mapsto\int\limits_{B(a)}u(x,t)\cdot w(x) dx
\end{equation}
strongly in $C([-S,0])$ for any $w\in L_2(B(a))$.
The corresponding sequences $v^{(k)}$ and $w^{(k)}$ converge  to their limits $v$ and $w$ in the same sense and of course $u=v+w$. For the pressure $p$, we have the following convergence: for any $n\in \mathbb N$,
there exists a sequences $c^{(k)}_n\in L_\frac 32(-S,0)$ such that
\begin{equation}\label{pressconv}
 \widetilde{p}^{(k)}_n\equiv   p^{(k)}-c^{(k)}_n\rightharpoonup p
\end{equation}
in $L_\frac 32(-S,0;L_\frac 32(B(n)))$.

So, arguing in the same way as in \cite{KS}, one can show that $u$ and  $p$ satisfy  the following conditions:
\begin{equation}\label{l2}
    \sup\limits_{-S<s<0}\sup\limits_{x_0\in\mathbb R^3}\|u(\cdot,s)\|^2_{L_2(B(x_0,1))}+
\sup\limits_{x_0\in\mathbb R^3}\int\limits_{-S}^0\int\limits_{B(x_0,1)}|\nabla u(y,s)|^2dyds<\infty;
\end{equation}
\begin{equation}\label{l3'}
    p\in L_\frac 32(-S,0;L_{\frac 32,{\rm loc}}(\mathbb R^3);
\end{equation}

the function
\begin{equation}\label{l4}
    s\mapsto \int\limits_{\mathbb R^3}u(y,s)\cdot w(y)dy
\end{equation}
is continuous on $[-S,0]$ for any compactly supported $w\in L_2(\mathbb R^3)$;
\begin{equation}\label{l5}
    \partial_tu+u\cdot\nabla u-\Delta u=-\nabla p,\quad {\rm div}\,u=0
\end{equation}
in $\mathbb R^3\times ]-S,0[$ in the sense of distributions;

for any $x_0\in \mathbb R^2$, there exists a function $c_{x_0}\in L_\frac 32(-S,0)$
such that
\begin{equation}\label{l6}
    p(x,t)-c_{x_0}(t)=p^1_{x_0}(x,t)+p^2_{x_0}(x,t)
\end{equation}
for all $x\in B(x_0,3/2)$;

for any $s\in ]-S,0[$ and for $\varphi\in C^\infty_0(\mathbb R^3\times ]-S,S[)$,
\begin{eqnarray}\label{l7}
% \nonumber to remove numbering (before each equation)
  \int\limits_{\mathbb R^3}\varphi^2(y,s)|u(y,s)|^2dy+2\int\limits_{-S}^s \int\limits_{\mathbb R^3}\varphi^2|\nabla u|^2dyd\tau\leq \nonumber\\
  \leq \int\limits_{-S}^s \int\limits_{\mathbb R^3}\Big(|u|^2(\Delta\varphi^2+\partial \varphi^2)+u\cdot\nabla \varphi^2(|u|^2+2p)\Big)dyd\tau.
\end{eqnarray}

Passing to the limit in (\ref{f13}), we find
$$\sup\limits_{x_0\in\mathbb R^3}\|v(\cdot,s)\|^2_{L_2(B(x_0,1))}\leq c(M)(s+S)^\frac 15$$
for all $s\in [-S,0]$. And thus
$$v\to 0 \qquad {\rm in}\quad L_{2,{\rm loc}}(\mathbb R^3)$$
as $s\downarrow -S$. Then, taking  into account (\ref{l1}), we can conclude that
\begin{equation}\label{l8}
 u\to a_0 \qquad {\rm in}\quad L_{2,{\rm loc}}(\mathbb R^3).
\end{equation}
as $s\downarrow -S$.

 By definition accepted in \cite{KS}, the pair $u$ and $p$, satisfying (\ref{l2})--(\ref{l8}), is a local energy solution to the Cauchy problem for the Navier-Stokes equations in $\mathbb R^3\times ]-S,0[$ associated with the initial velocity $a_0$.

%\begin{equation}\label{lp4}
%    p^{(k)}-c^{(k)}=\widetilde{p}^{(k)}\to p
%\end{equation}
%weakly in $L_\frac 32(B(a)\times ]-S,0[)$ for some suitable sequence $c^{(k)}\in L_\frac 32(-S,0)$.

Now, our aim is to show that $u$ is not identically zero.  %solution to the Navier-Stokes equations
%and, moreover, it is the so-called local energy solution in the interval $]-S_1,0[$ with some $S_1\leq %S$. Let us  start with the first task.
Using the inverse scaling, we observe that  the following identity takes place:
$$\frac 1{a^2}\int\limits_{Q(a)}(|u^{(k)}|^3+|\widetilde{p}^{(k)}|^\frac 32)dy\,ds=\frac 1{(a\lambda_k)^2}\int\limits_{Q(z_T,a\lambda_k)}(|v|^3+|q-b^{(k)}|^\frac 32)dx\,dt$$
for all
$0<a<a_*=\inf\{1, \sqrt{S/10},\sqrt{T/10}\}$ and for all $\lambda_k\leq 1$. Here, $z_T=(0,T)$, $\widetilde{p}^{(k)}\equiv \widetilde{p}^{(k)}_2$, and $b^{(k)}(t)=\lambda_k^{-2}c^{(k)}_2(s)$. Since the pair $v$ and $q-b^{(k)}$ is a suitable weak solution to the Navier-Stokes equations in $Q(z_T,\lambda_ka_*)$,  we find
\begin{equation}\label{lp5}
 \frac 1{a^2}\int\limits_{Q(a)}(|u^{(k)}|^3+|\widetilde{p}^{(k)}|^\frac 32)dy\,ds>\varepsilon
\end{equation}
for all $0<a<a_*$ with a positive universal constant $\varepsilon$.

Now, %we follow arguments from \cite{S6}. Our first observation is that,
by  (\ref{lp1}) and (\ref{pressconv}),
\begin{equation}\label{lp6}
 \frac 1{a^2}\int\limits_{Q(a)}|u^{(k)}|^3dy\,ds\to\frac 1{a^2}\int\limits_{Q(a)}|u|^3dy\,ds
\end{equation}
for all $0<a<a_*$ and
\begin{equation}\label{lp7}
   \sup\limits_{k\in \mathbb N} \frac 1{a_*^2}\int\limits_{Q(a_*)}(|u^{(k)}|^3+|\widetilde{p}^{(k)}|^\frac 32)dy\,ds= M_1<\infty.
\end{equation}
To treat the pressure $\widetilde{p}^{(k)}$, we do the usual decomposition of it into two parts,
see similar arguments in \cite{S6}. The first one is completely controlled by the pressure while the second one is a harmonic function in $B(a_*)$ for all admissible $t$. In other words,
we
have
$$\widetilde{p}^{(k)}=p^{(k)}_1+p^{(k)}_2$$
where $p^{(k)}_1$ obeys the estimate
\begin{equation}\label{lp8}
    \|p^{(k)}_1(\cdot,s)\|_{\frac 32,B(a_*)}\leq c\|u^{(k)}(\cdot,s)\|^2_{3,B(a_*)}.
\end{equation}
For the harmonic counterpart of the pressure, we have
$$\sup\limits_{y\in B(a_*/2)}  | p^{(k)}_2(y,s)|^\frac 32\leq c(a_*)\int\limits_{B(a_*)} | p^{(k)}_2(y,s)|^\frac 32dy$$
\begin{equation}\label{lp9}
\leq c(a_*) \int\limits_{B(a_*)}( |\widetilde{p}^{(k)}(y,s)|^\frac 32+|u^{(k)}(y,s)|^3)dy
\end{equation}
for all $-a_*^2<s<0$.

For any $0<a<a_*/2$,
$$\varepsilon\leq \frac 1{a^2}\int\limits_{Q(a)}( |\widetilde{p}^{(k)}|^\frac 32+|u^{(k)}|^3)dy\,ds\leq$$
$$\leq c \frac 1{a^2}\int\limits_{Q(a)}( |p_1^{(k)}|^\frac 32+|p_2^{(k)}|^\frac 32+|u^{(k)}|^3)dy\,ds\leq$$
$$\leq c \frac 1{a^2}\int\limits_{Q(a)}( |p_1^{(k)}|^\frac 32+|u^{(k)}|^3)dy\,ds+$$
$$+ca^3\frac 1{a^2}\int\limits_{-a^2}^0\sup\limits_{y\in B(a_*/2)}  | p^{(k)}_2(y,s)|^\frac 32ds.$$
From (\ref{lp7})--(\ref{lp9}), it follows that
$$\varepsilon\leq c\frac 1{a^2}\int\limits_{Q(a_*)}|u^{(k)}|^3dy\,ds+c a \int\limits_{-a^2}^0ds\int\limits_{B(a_*)}( |\widetilde{p}^{(k)}(y,s)|^\frac 32+|u^{(k)}(y,s)|^3)dy\leq$$
$$\leq c\frac 1{a^2}\int\limits_{Q(a_*)}|u^{(k)}|^3dy\,ds+ca\int\limits_{Q(a_*)}( |\widetilde{p}^{(k)}|^\frac 32+|u^{(k)}|^3)dy\,ds\leq$$
$$\leq c\frac 1{a^2}\int\limits_{Q(a_*)}|u^{(k)}|^3dy\,ds+cM_1aa^2_*$$
for all $0<a<a_*/2$. After passing to the limit and picking up sufficiently small $a$, we find
\begin{equation}\label{lp10}
    0<c\varepsilon a^2\leq \int\limits_{Q(a_*)}|u|^3dy\,ds
\end{equation}
for some positive $0<a<a_*/2$. So, the limit function $u$ is non-trivial.

\textsc{Proof Theorem \ref{mainresult}} %By (\ref{sb4}) and by the known compactness imbedding,
%we have
%$$u^{(k)}(\cdot,-S) \to u(\cdot,-S)$$
%in $L_{2,{\rm loc}}$. Moreover,
Since the limit function $a_0\in L_3$,
$$\|a_0\|_{2, B(x_0,1)}\to 0$$
as $|x_0|\to\infty$. The latter, together with Theorem 1.4 from \cite{KS} and $\varepsilon$-regularity theory for the Navier-Stokes equations, gives required decay at infinity.
The last thing to be noticed is that the following important property holds true:
\begin{equation}\label{lp19}
    u(\cdot,0)=0.
\end{equation}
This follows from (\ref{sb2}) and (\ref{lp1}), see the last statement in (\ref{lp1}).
 %the fact that in any case the quantity is bounded $\|v(\cdot, t_k)\|_3$ in $k$ and by $\varepsilon$-regularity theory, $\|v(\cdot, t_k)\|_3$ is bounded as well. Then applying ,
%So, we get (\ref{lp19}), f
More details on the matter can be found in papers \cite{S5} and \cite{S6}.
According to backward uniqueness for the Navier-Stokes, $u(\cdot,s)=0$  for any $s\in ]-a_*^2,0[$, which contradicts (\ref{lp10}). So, $z_T$ is not a singular point. Theorem \ref{mainresult} is proved.

\noindent\textbf{Acknowledgement} The author was partially supported by the RFFI grant 11-01-00324-a.

%\setcounter{equation}{0}
%\section{Preliminaries}

%In this section, we are going to recall the notion of local energy solutions introduced  by G.-P. %Lemarie-Rieusset, see monograph \cite{LR1} and references there. Here, we shall keep notation from  %paper \cite{KS} where a version of original statements from \cite{LR1} has been discussed.

\noindent

%G. Seregin\\
%Center for Nonlinear PDE's,\\
%Mathematical Institute, University of Oxford,UK\\
%and\\
%Steklov Institute of Mathematics at St.Petersburg, \\
%St.Peterburg, Russia


\begin{thebibliography}{99}
%\bibitem {CKN}
%Caffarelli, L., Kohn, R.-V., Nirenberg, L., Partial regularity of
%suitable weak solutions of the Navier-Stokes equations, Comm. Pure
%Appl. Math., Vol. XXXV (1982), pp. 771--831.

\bibitem {ESS4}

Escauriaza, L., Seregin, G.,  ~\v Sver\'ak, V.,
$L_{3,\infty}$-Solutions to the Navier-Stokes Equations and
Backward Uniqueness, Russian Mathematical Surveys, 58(2003)2, pp.
211-250.

%\bibitem{Kang}

%Kang, K., Unbounded normal derivative for the Stokes system near
%boundary, Math. Ann. 331(2005), pp. 87–-109.
%\bibitem {L}
%Ladyzhenskaya, O. A.,  Mathematical problems of the dynamics of
%viscous incompressible fluids, 2nd edition, Nauka, Moscow 1970.
\bibitem{Kato}

Kato, T., Strong $L_p$-solutions of the NavierпїЅStokes equation in $R^m$, with applications to
weak solutions. Math. Z., 187 (1984), 471--480.

%\bibitem {LS}
%Ladyzhenskaya, O. A., Seregin, G. A., On partial regularity of
%suitable weak solutions to the three-dimensional Navier-Stokes
%equations, J.  math. fluid mech.,  1(1999), pp. 356-387.




\bibitem{KNSS}
Koch, G., Nadirashvili, N., Seregin, G., Sverak, V.,
Liouville theorems for the Navier-Stokes equations and applications,
Acta Mathematica, 203 (2009), 83--105.

\bibitem {LR1}

Lemarie-Rieusset, P. G., Recent developemnets in the Navier-Stokes
problem, Chapman\&Hall/CRC reseacrh notes in mathematics series,
431.


\bibitem {Le}
Leray, J., Sur le mouvement d'un liquide visqueux emplissant
l'espace, Acta Math. 63(1934), pp. 193--248.

%\bibitem {Li}
%Lin, F.-H., A new proof of the Caffarelly-Kohn-Nirenberg theorem,
%Comm. Pure Appl. Math., 51(1998), no.3, pp. 241--257.

%\bibitem {Ne}
%Neustupa, J., Partial regularity of weak solutions to the
%Navier-Stokes equations in the class $L^\infty(0,T;L^3(\Om)^3)$,
%J. math. fluid mech., 1(1999), pp. 309--325.

\bibitem{KS}
Kikuchi, N., Seregin, G.,Weak solutions to the Cauchy problem for
the Navier-Stokes
 equations satisfying the local energy inequality, AMS translations, Series 2,
 Volume 220, pp.  141-164.

% \bibitem{KNSS}
% Koch, G., Nadirashvili, N., Seregin, G.,  ~\v Sver\'ak, V., Liouville
 %theorems for Navier-Stokes equations and applications, to appear in
 %Acta Mathematica.

 \bibitem{RS}
 Rusin, W., Sverak, V., Miminimal initial data for potential Navier-Stokes singularities,
arXiv:0911.0500.
%\bibitem {S1}
%Seregin, G.A., Some estimates near the boundary for solutions to
%the non-stationary linearized Navier-Stokes equations, Zapiski
%Nauchn. Seminar. POMI, 271(2000), pp. 204--223.

%\bibitem {S1'}
% Seregin, G., Local regularity theory of the Navier-Stokes equations, Handbook of Mathematical Fluid Mechanics, Vol. 4, Edited by Friedlander, D. Serre,
%pp, 159-200

%\bibitem {S1''}
%Seregin, G., A note on local boundary regularity for the Stokes system, to appear in Zapiski Nauchn. Seminar. POMI.



%\bibitem {Sc}
%Scheffer, V., Partial regularity of solutions to the Navier-Stokes
%equations, Pacific J. Math., 66(1976), 535--552.

%\bibitem {Sc1}
%Scheffer, V., Hausdorff measure and the Navier-Stokes equations,
%Commun. Math. Phys., 55(1977), pp. 97--112.

%\bibitem {S2}
%Seregin, G. A. On the number of singular points of weak solutions
%to the Navier-Stokes equations, Comm. Pure Appl. Math., 54(2001),
%issue 8, pp. 1019-1028.

%\bibitem {S3}
%Seregin, G.A., Local regularity of suitable weak solutions to the
%Navier-Stokes equations near the boundary,  J. math. fluid mech.,
%4(2002), no.1,1--29.


 %\bibitem {S4}
% Seregin, G.A.,   Differentiability properties of  weak solutions
%to the Navier-Stokes equations,  Algebra and Analysis, 14(2002),
%No. 1, pp. 193-237.

\bibitem{S5}
Seregin, G.A., Navier-Stokes equations: almost $L_{3,\infty}$-cases,
Journal of mathematical fluid mechanics, 9(2007), pp. 34-43.

\bibitem {S6}
Seregin, G., A note on necessary conditions for blow-up of energy solutions to the Navier-Stokes equations, Progress in Nonlinear Differential Equations
and Their Applications, 2011 Springer Basel AG, Vol. 60, 631--645.


\bibitem {S7}
Seregin, G., Necessary conditions of potential blow up for the Navier-Stokes equations, Zapiski Nauchn.Seminar. POMI, 385(2010), 187-199.
%\bibitem{SS1}

%Seregin, G.,  \v Sver\'ak, V.,  The Navier-Stokes equations and
%backward uniqueness, Nonlinear Problems in Mathematical Physics
%II, In Honor of Professor O.A. Ladyzhenskaya,  International
%Mathematical Series II, 2002, pp. 359--370.

\bibitem{SS3}
Seregin, G.,  \v Sver\'ak, V., On Type I singularities of the
local axi-symmetric solutions of the Navier-Stokes equations,
Communications in PDE's, 34(2009), pp. 171-201.

%\bibitem{SZ}
%Seregin, G., Zajaczkowski, W.,
 % A sufficient condition of regularity for axially symmetric
%  solutions to the Navier-Stokes equations,  SIMA J.  Math. Anal.,
 % (39)2007, pp. 669--685.

% \bibitem{Se}
%Serrin, J.,  On the interior regularity of weak solutions of the
%Navier-Stokes equations, Arch. Ration. Mech. Anal., 9(1962), pp.
%187--195.


\end{thebibliography}
\end{document}